\documentclass[12pt]{article}
\usepackage{amssymb}
\usepackage{amsfonts}
\usepackage{graphicx}

\newcommand{\Hil}[0]{
\mathcal{H} 
}
\newcommand{\norm}[2]{
\left\| #2 \right\|_{#1}
}
\newcommand{\HS}[0]{
{\mathcal HS}
}
\newcommand{\CC}[0]{
\mathbb{C} 
}
\newcommand{\qed}[0]{
$ \hspace{\stretch{1}} \Box $
}
\newcommand{\clsp}[1]{
{ \overline{span} \left\{ #1 \right\} }
}
\newtheorem{theorem}{Theorem}[section]
\newtheorem{def.}{Definition}[section]
\newtheorem{prop.}[theorem]{Proposition}
\newtheorem{lem.}[theorem]{Lemma}
\newtheorem{cor.}[theorem]{Corollary}
\newtheorem{Bsp.}{Example}[section]
\newcommand{\Bsp}[1]{
\begin{Bsp.} \label{#1} \bf : \end{Bsp.}
}
\newenvironment{proof}{\noindent \bf Proof: \rm}{\qed

}
\begin{document}

\title{Hilbert-Schmidt Operators and Frames - Classification, Approximation by Multipliers and Algorithms}

\author{
Peter Balazs\thanks{Acoustics Research Institute,
Austrian Academy of Science, Vienna, Austria
(Peter.Balazs@oeaw.ac.at). Supported partly by the EU network HASSIP, HPRN-CT-2002-00285.}
}
\maketitle
\begin{abstract}
In this paper we deal with the connection of frames with the class of Hilbert Schmidt operators. First we give an easy criteria for operators being in this class using frames. It is the equivalent to the criteria using orthonormal bases. Then we construct Bessel sequences frames and Riesz bases for the class of Hilbert Schmidt operators using the tensor product of such sequences in the original Hilbert space. We investigate how the Hilbert Schmidt inner product of an arbitrary operator and such a rank one operator can be calculated in an efficient way. Finally we give an algorithm to find the best approximation, in the Hilbert Schmidt sense, of arbitrary matrices by operators called frame multipliers, which multiply the analysis coefficents by a fixed symbol followed by a synthesis.

\vspace{5mm} 
\noindent {\it Key words and phrases} : frames, Bessel sequences, discrete expansion, operators, matrix representation, Hilbert-Schmidt operators, Frobenius matrices, Riesz bases.
\vspace{3mm}\\ 
\noindent {\it 2000 AMS Mathematics Subject Classification} (primary:) 41A58, 47B10, 65F30
\end{abstract}

\section{Introduction} \label{intro-sec}

Signal processing has become very important in today's life. Without signal processing methods several modern
technologies would not be possible, like mobile phone, UMTS, xDSL or digital television. 
The mathematical background for today's signal processing applications are \em Gabor \em \cite{feistro1} and \em wavelet \em \cite{daubech1} theory. 
%
From practical experience it soon has become 
apparent that the concept of an orthonormal basis is not always useful.
 This led to the concept of frames, which was introduced by Duffin and Schaefer \cite{duffschaef1}. It was made popular by Daubechies \cite{daubech1}, and today it is one of the most important foundations of Gabor \cite{feistro1} and wavelet \cite{aliant1} and sampling theory \cite{aldrgroech1}. In signal processing applications frames have received more and more attention \cite{boelc1,vettkov1}. 

Models in physics \cite{aliant1} and other application areas, for example in sound vibration analysis \cite{kreizxxl1}, are mostly continuous models. A lot of problems there can be formulated as operator theory problems, for example in differential or integral equations. 
A interesting class of operators is the Hilbert Schmidt class \cite{schatt1}. This class forms a Hilbert space and so is similar to Euclidean space in the sense, that there is a way to decide if two elements are orthogonal. This gives some tools to work with it, as for example the orthogonal projection, which is a very important and easy tool for example for finding minimal solutions.
This class of operators has therefore found also a lot of applications in signal processing, refer for example to \cite{miller00,shpa94}.

To be able to work numerically the operators have to be discretized, using matrices instead. Application and algorithms always work with finite dimensional data. So the numerical efficiency of involved algorithms is a very interesting questions, too. 

In this paper we will look at a connection of the theory of frames with the one of Hilbert Schmidt operators. It starts with prelimnaries and notations in Section \ref{sec:prelnot0}. In Section \ref{sec:framclasshs0} we will look at a way to describe these class of operators using frames. This can be done excately as in the case of orthonormal bases: An operator $H : \Hil \rightarrow \Hil$ is Hilbert Schmidt if and 
only if
$ \sum \limits_k \norm{\Hil}{H f_k}^2 < \infty $
for one (and therefore for all) frame(s) $(f_k)$.
In Section \ref{sec:frahMSOp0} we construct Bessel sequences frames and Riesz bases for the class of Hilbert Schmidt operators using tensor products of such sequences in the original Hilbert space. We will also investigate how the Hilbert Schmidt inner product of an arbitrary operator and such a rank one operator can be calculated in an efficient way, i.e. the inner product $\left< T , g_k \otimes f_l \right>_{\HS}$. This will be used in Section \ref{sec:bestapprhsframe0} to find the best approximation, in the Hilbert Schmidt sense, of an arbitrary matrix by a frame multiplier given an arbitary frame. For the regular Gabor case this has been investigated in \cite{feiham1}. We will give an algorithm working for all frame sequences. This is an important tool for applications: With the Gabor multipliers an efficient way to implement time-variant filters have been investigated \cite{hlawatgabfilt1}. Any time-variant linear system can
be modelled by a matrix, slowly time-varying linear system have a good
chance to be close to Gabor multipliers. 
Other matrices might be more connected to a 'diagonalization' using other frames. For application irregular Gabor frames \cite{liuwang1,xxlphd1}, wavelet frames \cite{daubech1} or other analysis / synthesis systems like e.g. Gammatone filter banks (e.g. refer to \cite{hartm1}), which are mainly used for analysis based on the auditory system, come to mind.
By the best approximation given in this paper we offer a method and an algorithm useable for all systems, that form a frame, to search for the best representation of any matrix as a diagonal operator using these systems.
We give MATLAB algorithms for this approximation. Those algorithms can be downloaded from \verb_http://www.kfs.oeaw.ac.at/xxl/finiteframes/hilschmidtfram.zip_.

\section{Preliminaries and Notations} \label{sec:prelnot0}

\subsection{Operators On Hilbert Spaces}
We will give only a short review, for details refer to \cite{conw1}. We will denote infinite dimensional Hilbert spaces by $\Hil$ 
and their inner product with $<.,.>$, which is linear in the first coordinate. 
An example for a Hilbert space is the sequence space $l^2$ consisting of all square-summable sequences in $\CC$
with the inner product $\left< c , d \right> = \sum \limits_{k} c_k \cdot \overline{d_k}$. 
%
%
Let $f \in \Hil_1$, $g \in \Hil_2$ then define the \em 
tensor product \em as an operator from $\Hil_2$ to $\Hil_1$ by
$ \left( f \otimes \overline{g} \right) (h) = \left< h, g \right> f  \mbox{ for } h \in \Hil_2$.

\subsubsection{Hilbert Schmidt operators} \label{sec:hilbertsch}

A bounded operator $T$ 
is called a \em Hilbert-Schmidt \em ($\HS$) operator if there exists an orthonormal basis (ONB) $( e_n ) \subseteq \Hil_1$ such that \index{operator!Hilbert Schmidt} 
$$ \norm{\HS}{T} := \sqrt{ \sum \limits_{n=1}^{\infty} \norm{\Hil_2}{T e_n}^2} < \infty .$$
Let ${\mathcal HS}(\Hil_1, \Hil_2)$ denote the space of Hilbert Schmidt operators from $\Hil_1$ to $\Hil_2$.
This definition is independent of the choice of the ONB.  The class of Hilbert-Schmidt operators is a Hilbert space of the compact operators with the following properties: 
\begin{itemize}
\item $\HS (\Hil_1, \Hil_2)$ is a Hilbert space with the inner product $\left< S , T \right>_{\HS} = \sum \limits_k \left< T e_k, S e_k \right>_\Hil$, where $( e_k)$ is an arbitrary ONB.
\item $\left< f \otimes \overline{g} , h \otimes \overline{l} \right>_\HS = \left< f, h\right>_\Hil \left< l, g \right>_\Hil$
\item $\norm{\HS}{T} = \norm{\HS}{T^*}$, and $T \in \HS$ $\Longleftrightarrow$ $T^* \in \HS$.
\end{itemize}

Furthermore it is easy \cite{xxlphd1} to show
\begin{cor.} \label{sec:HSinnprodmul1} Let $A \in \HS$, then
$\left< A , f \otimes \overline{g} \right>_\HS = \left< A g , f \right>_\Hil$.
\end{cor.}

For more details on this class of compact operators refer to \cite{schatt1} or \cite{wern1}.

\subsubsection{Frobenius Matrices} 

Let $\mathcal M_{m,n}$ be the vector space of $m \times n$-matrices with the inner product, 
$$\left< A , B \right>_{fro} = \sum \limits_{i=0}^{m-1} \sum \limits_{j=0}^{n-1} A_{i,j} \overline{B}_{i,j} \mbox{.}$$
called the space of {\em Frobenius matrices}. Extending this definition to infinite dimensional matrices, we call the space of matrices, where the corresponding norm converges, the space of Frobenius matrices, and will denote it by $l^{(2,2)}$ with the norm $\norm{2,2}{.}$.

Furthermore for $A$ be a $p \times q$ and $B$ a $r \times s$ matrix let the {\em Kronecker product} of $A$ and $B$ be the $p \cdot r \times q  \cdot s$ matrix $C$ with 
$$ C_{i,j} = a_{\left\lfloor \frac{i}{r} \right\rfloor,\left\lfloor \frac{j}{s} \right\rfloor} \cdot b_{i \  \rm mod \  \it r,j \  \rm mod \  \it s} $$
 
$$ A \otimes B = \left( \begin{array}{c c c c} a_{0,0} B & a_{1,0} B & \dots & a_{p-1,0} B \\
a_{0,1} B & a_{1,1} B & \dots & a_{n-1,1} B \\
\vdots & \vdots & \vdots & \vdots \\
a_{0,q-1} B & a_{1,q-1} B & \dots & a_{p-1,q-1} B
\end{array}  \right) $$

\subsection{Frames} \label{sec:fram0}
 
The sequence $\left( g_k | k \in K \right)$  is called a {\em frame} for the (separable) Hilbert space $\Hil$, if constants $A,B > 0$ exist, such that 
\begin{equation} \label{sec:framprop1} A \cdot \norm{\Hil}{f}^2 \le \sum \limits_k \left| \left< f, g_k \right> \right|^2 \le B \cdot  \norm{\Hil}{f}^2  \ \forall \ f \in \Hil
\end{equation}
$A$ is called a {\em lower} , $B$ a {\em upper frame bound}. 
The index set $K$ will be omitted in the following, if no distinction is necessary. 

If a sequence $(g_k)$ fulfills only the "upper frame condition", i.e.
$\sum \limits_k \left| \left< f, g_k \right> \right|^2 \le B \cdot  \norm{\Hil}{f}^2  \ \forall \ f \in \Hil $ 
it is called {\em Bessel sequence}.

%

Let $S_{\mathcal G} : \Hil  \rightarrow \Hil $ be the {\em frame  operator} 
$ S_{\mathcal G} ( f  ) = \sum \limits_k  \left< f , g_k \right> \cdot g_k$. 
%
%
$S = C^*C = DD^*$ is a positive invertible operator satisfying $A I_\Hil \le S \le B I_\Hil$ and $B^{-1} I_\Hil \le S^{-1} \le A^{-1} I_\Hil$.  
%
If we have a frame  in $\Hil$, we can find an expansion of every member of $\Hil$ with this frame.
Let $\left( g_k \right)$ be a frame for $\Hil$ with frame bounds $A$, $B > 0$. Then $\left( \tilde{g}_k \right) = \left( S^{-1} g_k \right)$  is a frame with frame bounds $B^{-1}$, $A^{-1} > 0$, the so called {\em canonical dual frame}. Every $f \in \Hil$ has an 
expansions
$ f = \sum \limits_{k \in K} \left< f, S^{-1} g_k \right> g_k $
and 
$ f = \sum \limits_{k \in K} \left< f, g_k \right> S^{-1} g_k $
where both sums converge unconditionally in $\Hil$.

Let $\{ g_k \}$ and $\{ f_k\}$ be two sequences in $\Hil$. 
The {\em cross-Gram matrix} $G_{g_k , f_k}$ for these sequences is given by $\left( G_{g_k, f_k} \right)_{jm} = \left< f_m , g_j \right>$, $j,m \in K$.  
If $(g_k) = (f_k)$ we call this matrix the {\em Gram matrix} $G_{g_k}$.

Using the pseudoinverse it can be shown \cite{olepinv} that  
$ D_{f_k}^\dagger = C_{\tilde{f_k}}$
and $ C_{f_k}^\dagger = D_{\tilde{f_k}}$. Furthermore it can be easily shown \cite{xxlphd1} that
\begin{cor.} Let $(g_k)$ be a frame on $\Hil$.
$$D_{g_k}^\dagger = C_{g_k} S_{g_k}^{-1} = G_{\tilde{g}_k,g_k} C_{\tilde{g}_k} = G_{g_k}^\dagger C_{g_k}$$
\end{cor.}

A complete sequence $( \psi_k)$ in $\Hil$  is called a \em Riesz basis \em if 
there exist constants $A$, $B >0$ such that the inequalities
$$ A \norm{2}{c}^2 \le \norm{\Hil}{\sum \limits_{k \in K} c_k \psi_k}^2 \le B \norm{2}{c}^2 $$
hold for all finite sequences $(c_k)$.
%
A frame $( \psi_k )$ is a Riesz basis for $\Hil$ if and only if 
there is a biorthogonal sequence $( \tilde{\psi}_k )$, i.e. $\left< \psi_k, \phi_j\right> = \delta_{kj}$ for all $h,j$. This is the dual frame.

For more details  on this topic and proofs refer to \cite{duffschaef1}, \cite{ole1} or \cite{Casaz1}.

\section{Frames Classifiying Hilbert-Schmidt Operators} \label{sec:framclasshs0}

It can be shown very easily that
\begin{lem.} Let $(f_k)$ be a frame and $(e_i)$ an ONB in $\Hil$. Let $H$ be an operator  $\Hil \rightarrow \Hil$. Then
$$ A \cdot \sum \limits_i \norm{\Hil}{H^* e_i}^2 \le \sum \limits_k \norm{\Hil}{H f_k}^2 \le B \sum \limits_i \norm{\Hil}{H^* e_i}^2 $$
\end{lem.}
\begin{proof}
$$\sum_k \norm{}{H f_k}^2 = \sum_k \sum_l \left| \left< H f_k, e_l\right> \right|^2 =  \sum_l \sum_k \left| \left< f_k, H^*e_l\right> \right|^2 = (*) $$
$$ (*) \ge A \sum_l \norm{}{H^* e_l}^2 $$
and 
$$ (*) \le B \sum_l \norm{}{H^* e_l}^2 $$
\end{proof}

From the proof it is clear that the right inequality is true for Bessel sequences, so
\begin{cor.} \label{sec:hilbschmibessel1} Let $(f_k)$ be a Bessel sequence and $(e_i)$ an ONB in $\Hil$. Let $H$ be an operator  $\Hil \rightarrow \Hil$. Then
$$ \sum \limits_k \norm{\Hil}{H f_k}^2 \le B \sum \limits_i \norm{\Hil}{H^* e_i}^2 $$
\end{cor.}
\vspace{5mm}

As we know that an operator is Hilbert Schmidt if and only if it's adjoint operator is as well, we get\footnote{Thanks and credits go to Ole Christensen. In private conversation he mentioned the existence of the first part of this corollary.}:
\begin{prop.}  \label{sec:framehs3}   An operator $H : \Hil \rightarrow \Hil$ is Hilbert Schmidt if and 
only if
$$ \sum \limits_k \norm{\Hil}{H f_k}^2 < \infty $$
for one (and therefore for all) frame(s).

$$ \sqrt{A} \norm{\HS}{H} \le \sqrt{\sum \limits_k \norm{\Hil}{H f_k}^2} \le \sqrt{B} \norm{\HS}{H} $$ 

In particular for tight frames we have 
$$ \norm{\HS}{H} = \frac{1}{A} \sqrt{\sum \limits_k \norm{\Hil}{H f_k}^2}$$ 
\end{prop.}

\section{Frames In The Hilbert-Schmidt Class Of Operators} \label{sec:frahMSOp0}

For two frames $\{ f_k\}, \{g_k\}$ an operator can be described \cite{xxlframoper1} by the matrix 
$$ \mathcal{M}^{(f_k , \tilde{g}_j)}  \left( O \right)_{k,j} = 
\left<O \tilde{g}_j, f_k \right>.$$
 This is the matrix that maps $C_{f_k} ( f ) \mapsto C_{ g_k} ( f )$. It is identical to the $\HS$ inner product of $O$ and $f_k \otimes \overline{\tilde{g}_j}$.
 This can be applied to Hilbert-Schmidt operators \cite{xxlframoper1} as for $O \in \HS$ the representing $\mathcal{M}(O)$ is a Frobenius matrix
 with $\norm{2,2}{\mathcal{M}(O)} \le \sqrt{B B'} \norm{\HS}{O}$. This mapping $\mathcal{M}^{(f_k , \tilde{g}_j)} : \HS \rightarrow l^{(2,2)}$ is injective.

Also the opposite mapping can be investigated. For a Frobenius matrix $M$ define 
$\left( \mathcal{O}^{(\Phi , \Psi)} \left( M \right)\right) h = \sum \limits_k  \left( \sum \limits_j M_{k,j} \left<h, \psi_j\right> \right) \phi_k \mbox{, for } h \in \Hil_1$
then $ \mathcal{O}(M)$ is in $\HS$ and $\norm{\HS}{ \mathcal{O}(M) } \le \sqrt{B B'} \norm{fro}{M}$, cf. \cite{xxlframoper1}.

If the involved frames are Riesz bases these two operators are bijective.
This can be used to investigate frames in the class of Hilbert-Schmidt operators:

\begin{theorem}  \label{sec:hsframeriesz1} Let $(g_k)_{k \in K}$ be a sequence in $\Hil_1$ , $(f_i)_{i \in I}$ in $\Hil_2$. Then
\begin{enumerate}
\item Let $(g_k)$ and $(f_i)$ be Bessel sequences with bounds $B$,$B'$, then $( f_i \otimes \overline{g}_k )_{(i,k) \in I \times K}$ is a Bessel sequence for $\HS(\Hil_1,\Hil_2)$ with bound $\sqrt{B \cdot B'}$.
\item Let $(g_k)$ and $(f_i)$ be frames with bounds $A,B$ and $A,B'$. Then $( f_i \otimes \overline{g}_k )_{(i,k) \in I \times K}$ is a frame for $\HS(\Hil_1,\Hil_2)$ with bounds $\sqrt{A \cdot A'}$ and $\sqrt{B \cdot B}$. A dual frame is $( \tilde{f}_i \otimes \overline{\tilde{g}}_k )$.
\item Let $(g_k)$ and $(f_i)$ be Riesz bases. Then $( f_i \otimes \overline{g}_k )_{(i,k) \in I \times K}$ is a Riesz basis for $\HS(\Hil_1,\Hil_2)$. The biorthogonal sequence is $( \tilde{f}_i \otimes \overline{\tilde{g}}_k )$.
\end{enumerate}
\end{theorem} 
\begin{proof}
Suppose the operator $O \in \HS$, then 
$$ \mathcal{M}^{(f_k,g_j)} \left( O \right)_{k,j} = \left<O g_j, f_k \right>_{\Hil_1} \stackrel{Cor. \ref{sec:HSinnprodmul1}}{=} \left<O , f_k \otimes \overline{g}_j \right>_{\HS} $$
Using the dual frames we know \cite{xxlframoper1} that 
$$ \norm{\HS}{O} = 
\norm{\HS}{\mathcal{O}^{(\tilde{f}_k,\tilde{g}_j)} \left(  \mathcal{M}^{(f_k,g_j)} \left(O\right) \right) } \le 
\frac{1}{\sqrt{A A'}} \norm{\HS}{M^{(f_k,g_j)}(O)},$$
as the reciprocal of the lower frame bound is an upper frame bound of the dual frame.
Therefore
$$ A A' \norm{\HS}{O}^2 \le \norm{\HS}{M^{(f_k,g_j)}(O)}^2 \le B B' \norm{\HS}{O}^2 $$
This is equal to
$$ A A' \norm{\HS}{O}^2 \le \sum \limits_{k,j} \left| \left<O , f_k \otimes \overline{g}_j \right>_{\HS} \right|^2 \le B B' \norm{\HS}{O}^2 $$

%
If both sequences are Riesz bases, then
$$ \left< f_{k_1} \otimes \overline{g}_{j_1} , \tilde{f}_{k_2} \otimes \overline{\tilde{g}}_{j_2} \right>_{\HS} = \left< f_{k_1} , \tilde{f}_{k_2} \right>_{\Hil} \cdot  \left< g_{j_2} , g_{j_1} \right>_{\Hil} = $$
$$ = \delta_{k_1,k_2} \cdot \delta_{j_1,j_2}. $$
This means that $\left( g_k \otimes \overline{f}_j \right)$ is a Riesz basis.
\end{proof}


\subsection{Efficient Calculation of the $\HS$ Inner Product}

Let us look at the finite dimensional space $\CC^n$. For matrix representation of operators \cite{xxlframoper1} in $\HS$, the inner product $\left< T , g_k \otimes f_l \right>_{\HS}$ becomes important. The diagonal version $\left< T , g_k \otimes f_k \right>_{\HS}$ plays an essential role for frame multipliers, see Section \ref{sec:bestapprhsframe0}. There are several ways to calculate this $\HS$ inner product, which we will list in Theorem \ref{sec:calcHSinnerprod1}. We will first collect the following properties for the proof of this theorem.

Note that with $\left\lfloor x \right\rfloor$ we describe the biggest integer smaller than $x$.\index{symbols!$\left\lfloor x \right\rfloor$}

\begin{lem.} \label{sec:matr2vec1} 
For $M \in \mathcal M_{m,n}$, a Frobenius matrix of dimension $m \times n$, let
$$ \mathfrak{vec}^{(n)} \left( M \right)_{k} = M_{k \ \rm mod \ \it n , \left\lfloor \frac{k}{n} \right\rfloor} \ for \ k = 0,\dots,m\cdot n - 1$$ \index{symbols!${\mathfrak vec}^{(n)}{M}$}
With this function this space is isomorphic to $\CC^{m \times n}$ with the standard inner product. The inverse of this function is
$$ \mathfrak{Mat_n} (x)_{i,j} = x_{i + j \cdot n}$$ \index{symbols!$\mathfrak{Mat_n} (x)$}
\end{lem.}
\begin{proof} The function $\mathfrak{vec}^{(n)}$ is clearly linear and inverse to $\mathfrak{Mat_n}$.
\end{proof}
The function $\mathfrak{vec}^{(n)}$ joins the columns together to a vector. The function $x \mapsto \mathfrak{Mat_n} (x)$ separates this vector again. $\mathfrak{Mat_n} (x)$ is well known in signal processing, it is called \em Polyphase representation \em there. 
\\

We will use the standard notation for complexity of a formula $\mathcal{O}$\cite{knuth97}, 
writing $f(x) \in O(n)$ if the algorithm for $f(x)$ has order of $n$ time complexity. We will abuse the standard notation a bit, as we will not only include the highest power term.

\begin{lem.} \label{sec:complebasic1} The complexity of the calculation of the following terms is
\begin{enumerate}
\item inner product: Let $x,y \in \CC^p$, then
$$\left< x , y \right> \in \mathcal{O}\left( 2p-1 \right)$$
\item matrix-vector multiplication: Let $A \in \mathcal M_{p,q}$, $x \in \CC^q$, then
$$A \cdot x \in \mathcal{O}\left( p \cdot \left( 2 q -1 \right) \right)$$
\item matrix-matrix multiplication: Let $A \in \mathcal M_{p,q}$, $B \in \mathcal M_{q,r}$, then
$$A \cdot B \in \mathcal{O}\left( p \cdot r \cdot \left( 2 q -1 \right)\right)$$
\item Kronecker product of matrices: Let $A \in \mathcal M_{p,q}$, $B \in \mathcal M_{r,s}$, then
$$A \otimes B \in \mathcal{O}\left( p \cdot q \cdot r \cdot s \right)$$
\end{enumerate}
\end{lem.}
\begin{proof} This can be easily be deduced from the basic definitions of these operators.
\end{proof}

\begin{lem.} \label{ref:lemkronmatrix1} Let $A \in \mathcal M_{r,s}$, $B \in \mathcal M_{p,q}$ and $C \in \mathcal M_{q,r}$. Then
$$ \left( A^T \otimes B \right) \cdot \left( \mathfrak{vec}^{(q)} C \right) = \mathfrak{vec}^{(p)}\left(B \cdot C \cdot A\right)$$
\end{lem.}
\begin{proof}
$$ \left( \left( A^T \otimes B \right) \left( \mathfrak{vec}^{(q)} C \right) \right)_i = \sum \limits_{j = 0}^{q \cdot s -1} \left( A^T \otimes B \right)_{i,j} \left( \mathfrak{vec}^{(q)}{C}\right)_j  = $$
$$ \sum \limits_{j = 0}^{q \cdot s -1}  A^T_{\left\lfloor \frac{i}{p} \right\rfloor,\left\lfloor \frac{j}{q} \right\rfloor} \cdot  B_{i \ \rm mod \ \it p ,j \ \rm mod \ \it q} C_{j \ \rm mod \ \it q ,  \left\lfloor \frac{j}{q} \right\rfloor} = (*)$$
Let $j_1 = j \ \rm mod \ \it q$ and $j_2 = \left\lfloor \frac{j}{q} \right\rfloor$, so
$$ (*) = \sum \limits_{j_1 = 0}^{q - 1} \sum \limits_{j_2 = 0}^{s -1}  A^T_{\left\lfloor \frac{i}{p} \right\rfloor,j_2} \cdot  B_{i \ \rm mod \ \it p ,j_1 } C_{j_1,j_2} = $$
$$ = \sum \limits_{j_2 = 0}^{s -1}  A_{j_2,\left\lfloor \frac{i}{p} \right\rfloor} \cdot  \left(B \cdot C\right)_{i \ \rm mod \ \it p,j_2} = $$
$$ = \left( B \cdot C \cdot A \right)_{i \ \rm mod \ \it p , \left\lfloor \frac{i}{p} \right\rfloor}$$
\end{proof}

\begin{theorem} \label{sec:calcHSinnerprod1}
Let $(h_l)_{l=0}^{L}$ be a frame in $\CC^n$, $(g_k)_{l=0}^{K}$ in $\CC^m$. Let $T$ be a linear operator $T : \CC^n \rightarrow \CC^m$. Then 
\begin{enumerate} 
\item $\left< T , g_k \otimes \overline{h}_l \right>_{\HS} = \left< \mathfrak{vec}^{(n)} \left( T \right) , \mathfrak{vec}^{(n)} \left( g_k \otimes \overline{f}_l \right) \right>_{\CC^{m \cdot n}} $ $ \in \mathcal{O}\left( (3 m n + m -1) \right)$ for each single pair $(l,k)$.
\item $\left< T , g_k \otimes h_l \right>_{\HS} = \left< T h_l , g_k \right>_{\CC^{m}} $ $\in \mathcal{O}\left( (2 m n + m -1) \right)$ for each single pair $(l,k)$.
\item $\left< T , g_k \otimes \overline{h}_l \right>_{\HS} = \left( C_{g_k} \cdot T \cdot D_{h_l} \right)_{l,k} $ $\in \mathcal{O}\left( \left( L \left( 2 mn - m + 2 m K -K \right) \right) \right)$ for all values $(l,k)$.
\item $\left< T , g_k \otimes h_l \right>_{\HS} = \left( D^T_{g_k}  \otimes C_{f_l} \right) \mathfrak{vec}(T) $ $\in \mathcal{O}\left( \left( KL \cdot \left( 3 m n - 1 \right) \right)) \right)$ for all values $(l,k)$.
\end{enumerate}
\end{theorem}
\begin{proof} We will use Lemma \ref{sec:complebasic1} extensively:

1.) $g_k \otimes h_l$ $\in \mathcal{O}\left( m \cdot n \right)$! $\mathfrak{vec}^{(n)}$ is only a reordering. The complex conjugation $\in \mathcal{O}\left( m \right)$. 
The inner product $\in \mathcal{O}\left( 2\cdot (m n) -1 \right) $. So the overall sum is $3 mn + m -1$.

2.) $ T h_l$ $\in \mathcal{O}\left( m \left( 2 n -1\right) \right)$. The inner product $\in \mathcal{O}\left( 2 m -1 \right)$. The sum is  $2 m n +m -1$.
 
3.) To understand that $\left< T , g_k \otimes h_l \right>_{\HS} = \left( C_{g_l} \cdot T \cdot D_{h_l} \right)_{l,k} $ use Corollary \ref{sec:HSinnprodmul1} and the facts that
$$ T \cdot \left( \begin{array}{c c c c } \vline & \vline & & \vline \\
g_1 & g_2 & \dots & g_M  \\
\vline & \vline & & \vline \\
 \end{array} \right) = \left( \begin{array}{c c c c } \vline & \vline & & \vline \\
T g_1 & T g_2 & \dots & T g_M  \\
\vline & \vline & & \vline \\
 \end{array} \right) $$
and
$$ \left( \left( \begin{array}{c c c}\mbox{---} & \overline{h_1} & \mbox{---} \\
\mbox{---} & \overline{h_2} & \mbox{---} \\
\vdots & & \vdots \\
\mbox{---} & \overline{h_N} & \mbox{---} \\
\end{array} \right) \cdot 
 \left( \begin{array}{c c c c } \vline & \vline & & \vline \\
g_1 & g_2 & \dots & g_M  \\
\vline & \vline & & \vline \\
 \end{array} \right)  \right)_{m,n} = \left< g_n , h_m \right> $$
Furthermore
$T \cdot D \in \mathcal{O}\left( m L (2 n - 1)\right)$, $C \cdot \left( T D \right) \in \mathcal{O}\left( K L (2 m -1) \right)$, so altogether we get as a sum $  m L 2 n - m L + K L 2 m - K L $ $= L \left( 2 mn - m + 2 m K -K\right)$. 

4.) Using Lemma \ref{ref:lemkronmatrix1} we know that this equality is true. For the calculation of $D^T_{g_k}  \otimes C_{f_l}$ we need 
$\in \mathcal{O}\left( K m L n \right)$. And for the matrix vector multiplication in $\CC^{mn}$ $\in \mathcal{O}\left( K L (2 m n -1 )\right)$. So overall $K L (2 m n -1 ) + K m L n = KL \cdot \left( 3 m n - 1 \right)$.\end{proof}

So overall if we have to calculate the inner products for all pairs $(k,l)$ the third method is the fastest (except when $n$ is very big and $m$ and $K$ very small). If we need only the diagonal part $k = l$, the second one is the most efficient as for using the third method we would still have to calculate the whole matrix and then use its trace.

\section{Best Approximation of Hilbert Schmidt Operators by Frame Multipliers} \label{sec:bestapprhsframe0}

\subsection{Best Approximation By Frame Sequences} \label{sec:apprmatrframmul0}

In applications we very often have the problem to find an approximation of a certain object. 
In a Hilbert space setting, where the interesting objects are in a space spanned by a frame sequence, it is known \cite{ole1} that the best approximation is just the orthogonal projection on this space, which is in this case given by
$$ P(f) = \sum \limits_k \left< f , \tilde{g}_k \right> g_k  = (*) $$ 

A disadvantage of this formula for practical solutions is that the dual frame has to be calculated. This can be time-consuming and is not needed per-se. But we can use the formulas we have established previously to get
$$ (*) = D_{g_k} C_{\tilde{g}_k} =  D_{g_k} D_{\tilde{g}_k}^\dagger = D_{g_k} G_{g_k}^\dagger C_{g_k} $$
Here we can avoid calculating the dual frame directly, instead using the existing algorithm to calculate the pseudoinverse. 

\begin{theorem} \label{sec:bestapprframseq1} Let $(g_k)$ be a frame sequence in $\Hil$. Let $V = \overline{span( g_k )}$. The best approximation of an arbitrary element $f \in \Hil$ is
$$ P(f) = D_{g_k} G_{g_k}^\dagger C_{g_k} f $$
\end{theorem}


\subsection{Frame Multiplier}

Frame Multipliers \cite{xxlmult1} are operators which can be described in this sense by diagonal matrices, i.e. operators of the form
$$ {\bf M}_{m, f_k, g_k} =  \sum \limits_k m_k \left< f, g_k \right> f_k,$$
where $(m_k)$ is a bounded sequence called the \em upper symbol\em . 
These are operators spanned by $\left( \gamma_k \otimes \overline{g}_k \right)$. We know from Theorem \ref{sec:hsframeriesz1} for the Hilbert-Schmidt class of operators 
\begin{enumerate}
\item that this is a Bessel sequence for Bessel sequences $(g_k)$ and $(f_k)$
\item and because every sub-family of a Riesz basis is a Riesz sequence, that for Riesz bases, $(\gamma_k \otimes \overline{g}_k )$ is a Riesz sequence.
\end{enumerate}

\subsubsection{The Identity As Multiplier} \label{sec:tightid}

Let us state an interesting question regarding multipliers, which can be adressed by the approximation algorithm presented later. Can the identity be described as multiplier?
It can be easily shown \cite{xxlphd1} that
\begin{lem.} \label{sec:tightid1} If and only if the identity is a multiplier for the Bessel sequence $\{ g_k \}$ with constant symbol $c \not= 0$, then $\{g_k\}$ is a tight frame.
\end{lem.}

In the case of regular well-balanced Gabor frames, it is shown in \cite{doerf1} that if the identity can be written as multiplier for the frame $\{g_k\}$, its symbol is a constant sequence. In this case Lemma \ref{sec:tightid1} is clearly equivalent to: The identity is a multiplier if and only if $\{g_k\}$ is a tight frame. We will see some examples later, that this is certainly not the case for general frames.

\subsection{The Lower Symbol} \label{sec:lowsymb0}

For Bessel sequences $(g_k)$ and $(f_k)$ the family $( g_k \otimes \overline{f}_k)$ is again a Bessel sequence.  So the synthesis operator is well defined:
$$ C_{g_k \otimes \overline{f}_k} : \HS \rightarrow l^2 \mbox{ with } C(T) = \left< T, g_k \otimes \overline{f}_k \right>_\HS $$
Using Theorem \ref{sec:calcHSinnerprod1} we can express that inner product as
$$C(T) = \left< T, g_k \otimes \overline{f}_k \right>_\HS = \left< T \overline{f}_k, g_k \right>_\Hil $$
So we define
\begin{def.} Let $( g_k )$ and $( \overline{f}_k )$ be Bessel sequences for $\Hil$, then the \bf lower symbol \it of an operator $T \in \HS$ is defined as 
$$ \sigma_L (T) = \left< T f_k , g_k \right>_\Hil $$ 
\end{def.}

The function $\sigma_L : \HS \rightarrow l^2 $ is just the synthesis operator of the Bessel sequence $g_k \otimes \overline{f}_k$ in $\HS$ and therefore well defined in $l^2$. The name is used in acoordance with the use for Gabor multipliers, which are frame multipliers for Gabor frames, refer for example to \cite{feinow1}. 

Let us suppose that 
these elements $(g_k \otimes \overline{f_k})$ are not only a Bessel sequence, but even a frame sequence. This has to be supposed and cannot not be deduced from the status of $(f_k)$ and $(g_k)$ of being frames, as Theorem \ref{sec:hsframeriesz1} is of no use here. Subsequences of frames don't have to be frames again \cite{ole1}.

In this case the best approximation can be found by using the analysis and the dual synthesis operator for the projection on the closed span of the elements $V = \clsp{ g_k \otimes \overline{f_k}}$, which are exactly those $\HS$ operators that can be expressed as frame multipliers with the given frames.
Let $Q_k$ be the canonical dual frame of $g_k \otimes \overline{f}_k$ in $V$ then the best approximation is
$$ P_V (T) = \sum \limits_k \left< T , g_k \otimes \overline{f}_k \right>_\HS Q_k = \sum \limits_k \sigma_L (T) Q_k $$
As $\left< T , g_k \otimes \overline{f}_k \right>_\HS$ is just the analysis coefficients with a frame sequence, we know $\norm{2}{\sigma_U} \le \norm{2}{\sigma_L (T)}$ for any other coefficients $\sigma_U$ such that the projection $P_V$ can be expressed in this way
and hence the name 'lower symbol'. Also for bounded operators $T$ which are not in $\HS$ this inner product is defined and bounded by $\norm{Op}{T} \sqrt{B B'}$. 

\subsubsection{Approximation Of Matrices By Frame Multipliers} \label{sec:apprframmult0}

In infinite-dimensional spaces not every subsequence of a frame is a frame sequence, but in the finite-dimensional case, all sequences are frame sequences. So we can use the ideas in the Sections \ref{sec:apprmatrframmul0} and \ref{sec:lowsymb0}  and apply it to frame multipliers. All codes in here can be found online at\\ \verb_http://www.kfs.oeaw.ac.at/xxl/finiteframes/hilschmidtfram.zip_. 

We want to find the best approximation (in the Frobenius norm) of a $m \times n$ matrix $T$ by a frame multiplier with the frames $(g_k)_{k=1}^K \subseteq \CC^n$ and $(f_k)_{k=1}^K \subseteq \CC^m$. This whole section is a generalization of the ideas in \cite{feiham1}.
\\

\verb_Algorithm:_

\begin{enumerate}
\item \verb_Inputs: T, D, Ds_ \\
$T$ is a $m \times n$ matrix, $D$ is the $n \times K$ synthesis matrix of the frame $(g_k)$. This means \cite{xxlfinfram1} that the elements of the frame are the columns of $D$. $Ds$ is the synthesis matrix of the frame $(f_k)$. Often we will use the case $(f_k) = (g_k)$ so $D_s = D$ by default.
\item \verb_Lower Symbol :_ \\
Using Lemma \ref{sec:complebasic1} the most efficient way to calculate the inner product $\left< T , g_k \otimes \overline{f}_k \right>_\HS$ is $\left< T f_k , g_k \right>_{\CC^n}$. This can be implemented effectively using the matrix multiplication by

\verb! (MATLAB :)     lowsym(i) = conj(D(:,i)'*(T*Ds(:,i)));!
\item \verb_Hilbert Schmidt Gram Matrix :_ \\
We calculate the Gram matrix of $(g_k \otimes \overline{f}_k)$ 
$$ {\left( G_{\HS} \right)}_{l,k} = \left< g_k \otimes \overline{f}_k , g_l \otimes \overline{f}_l \right>_\HS = \left< g_k , g_l \right>_\Hil \cdot \left< f_l , f_k \right>_\Hil = \left( G_{g_k} \right)_{l,k} \cdot \left( G_{f_k} \right)_{k,l}  $$
\verb! (MATLAB :) Gram = (D'*D).*((Ds'*Ds).');! \\
If $(g_k) = (f_k)$ then
$$ {\left( G_{\HS}\right)}_{l,k} = \left| \left< g_k , g_l \right>_\Hil \right|^2 $$
It is more efficient to use this formula in \\
\verb! (MATLAB :) Gram = abs((D'*D)).^2;! \\
as this has complexity, using Lemma \ref{sec:complebasic1}, $\sim K^2 \cdot \left( n^2 +2 \right)$ compared to the original calculation with $\sim K^2 \cdot \left( n^2 + m^2 +1\right)$.
\item \verb_Upper Symbol :_ \\
Using Theorem \ref{sec:bestapprframseq1} we get the coefficients of the approximation by using the pseudoinverse of the Gram matrix. In the case of frame multipliers the coefficients are an upper symbol $\sigma$. \\
\verb! (MATLAB :) uppsym = pinv(Gram)*lowsym;! \\
\item \verb_Outputs: TA, uppsym_ \\
For the calculation of the approximation we just have to apply the synthesis operator of the sequence $(g_k \otimes \overline{f}_k)$ to the upper symbol.
$$ \verb_TA_ = P_V(T) = \sum \limits_{k=1}^K \sigma_k g_k \otimes \overline{f}_k $$
The matrix of the operator $g_k \otimes \overline{f}_k$ can just be calculated \cite{xxlfinfram1} by ${\left( g_k \right)}_i \cdot  {\left( \overline{f}_k \right)}_j$. \\
\verb! (MATLAB :) P = D(:,i)*Ds(:,i)'; ! \\
\end{enumerate}
For an implementation of this algorithm in MATLAB see the file \verb_ApprFramMult.m_.

\Bsp{sec:testapprframmult} We will look at two simple example in $\CC^2$.
\begin{enumerate}
\item Let $A = \left( \begin{array}{c c} 3 & 0 \\ 0 & 5 \end{array} \right)$. This is clearly a multiplier for the standard orthonormal basis of $\CC^2$. The sequence $f_1 = \left( \frac{1}{2} , \frac{\sqrt{3}}{2} \right) , f_2 = \left( \frac{\sqrt{3}}{2} , - \frac{1}{2} \right)$ is also an ONB. But the best approximation of $A$ with this basis is $P_V (A) = \left( \begin{array}{c c} 3.7500 & 0.4330 \\  0.4330 &   4.2500 \end{array} \right)$. So this is an example that not even for ONBs a frame multiplier for one basis stays a frame multiplier for another one.
\item Let $T = Id_{\CC^2}$. and let $D = \left( \begin{array}{c c c} \cos (30^\circ) & 1 &  0 \\ \sin (30^\circ) &  1  &-1 \end{array} \right)$. This is a frame with bounds $A = 0.5453$, $B = 3.4547$ and therefore not tight. Still the identity can be approximated perfectly (up to numerical precision) with the coefficients $\sigma = \left( 3.1547 , -1.3660, 1.5774\right)$. So this is an example, where the identity is a frame multiplier for a non-tight system.
\end{enumerate}

The MATLAB-codes for these examples can be found in \verb_testapprframmult.m_.

\begin{figure}[!ht]
	\begin{center}
		\includegraphics[width=0.49\textwidth]{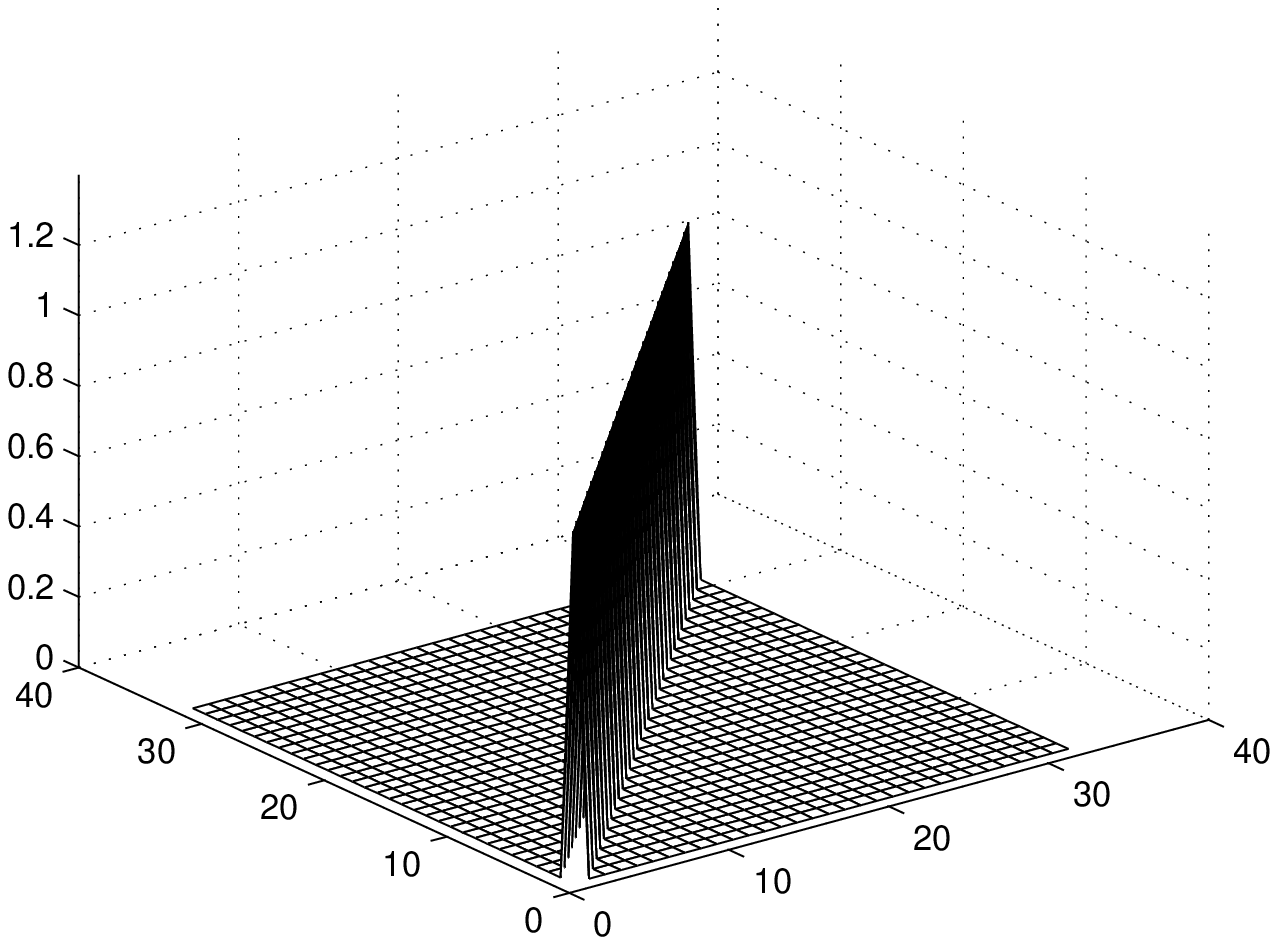}
		\includegraphics[width=0.49\textwidth]{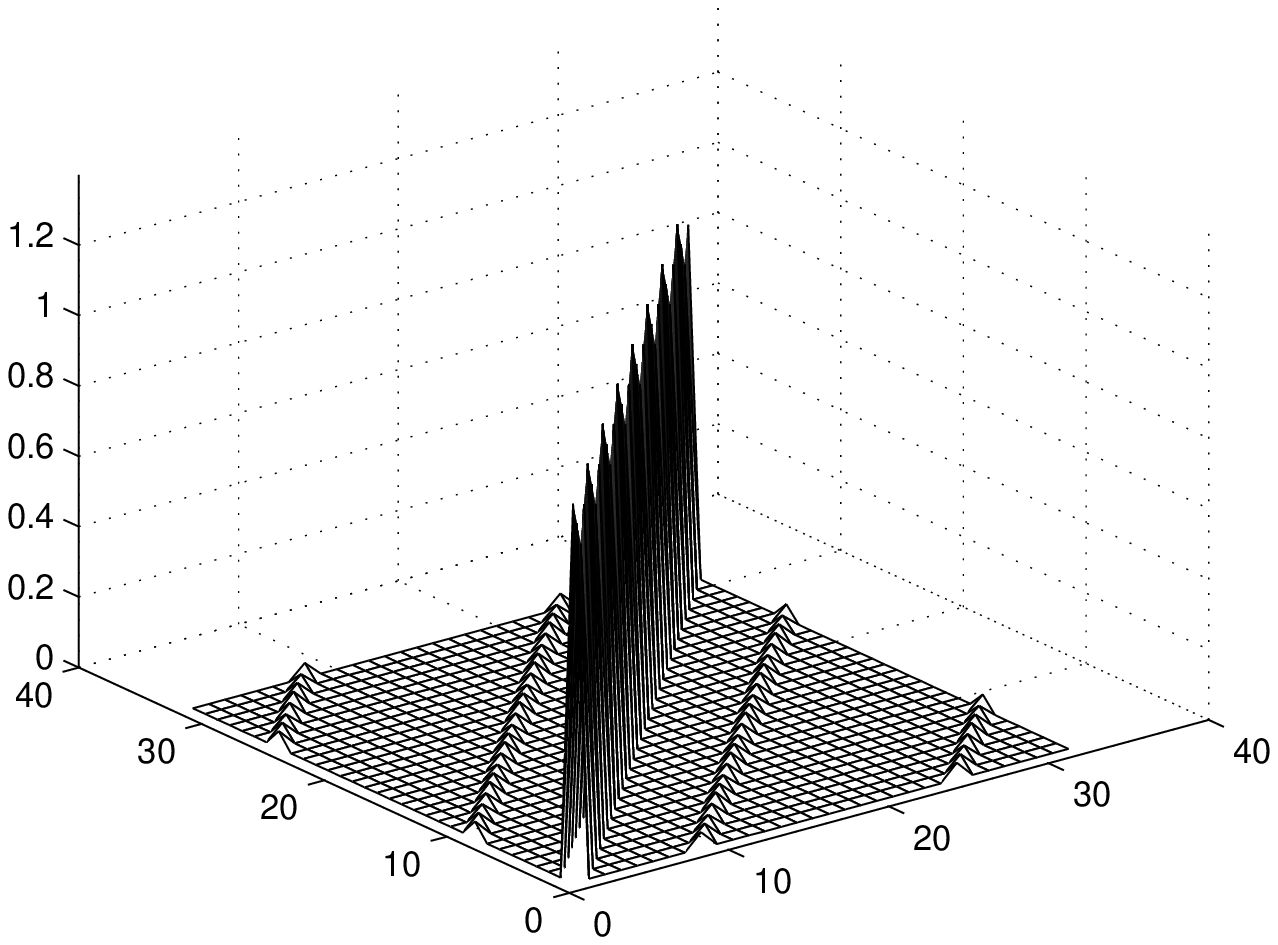}
		\includegraphics[width=0.49\textwidth]{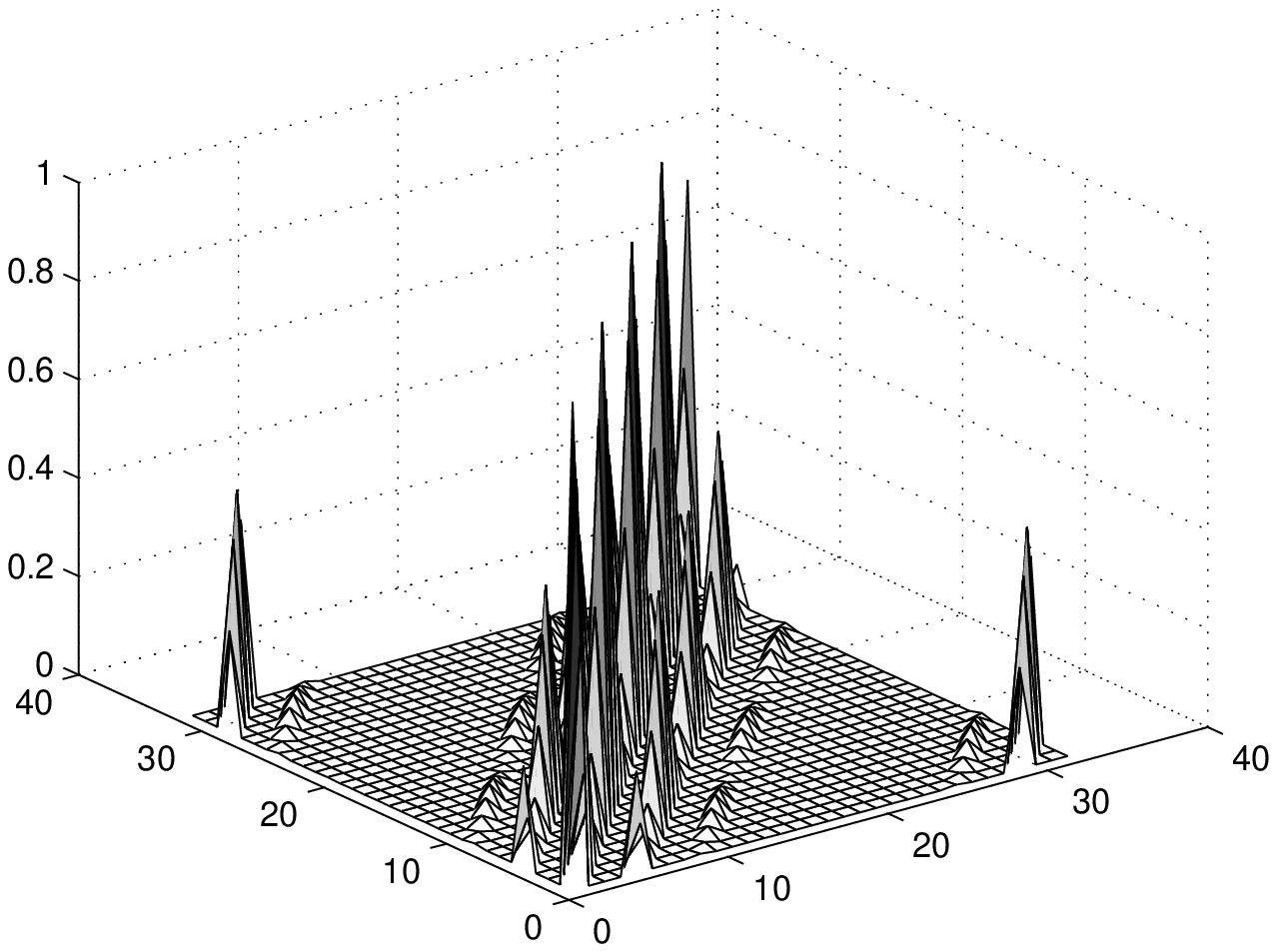}
		\includegraphics[width=0.49\textwidth]{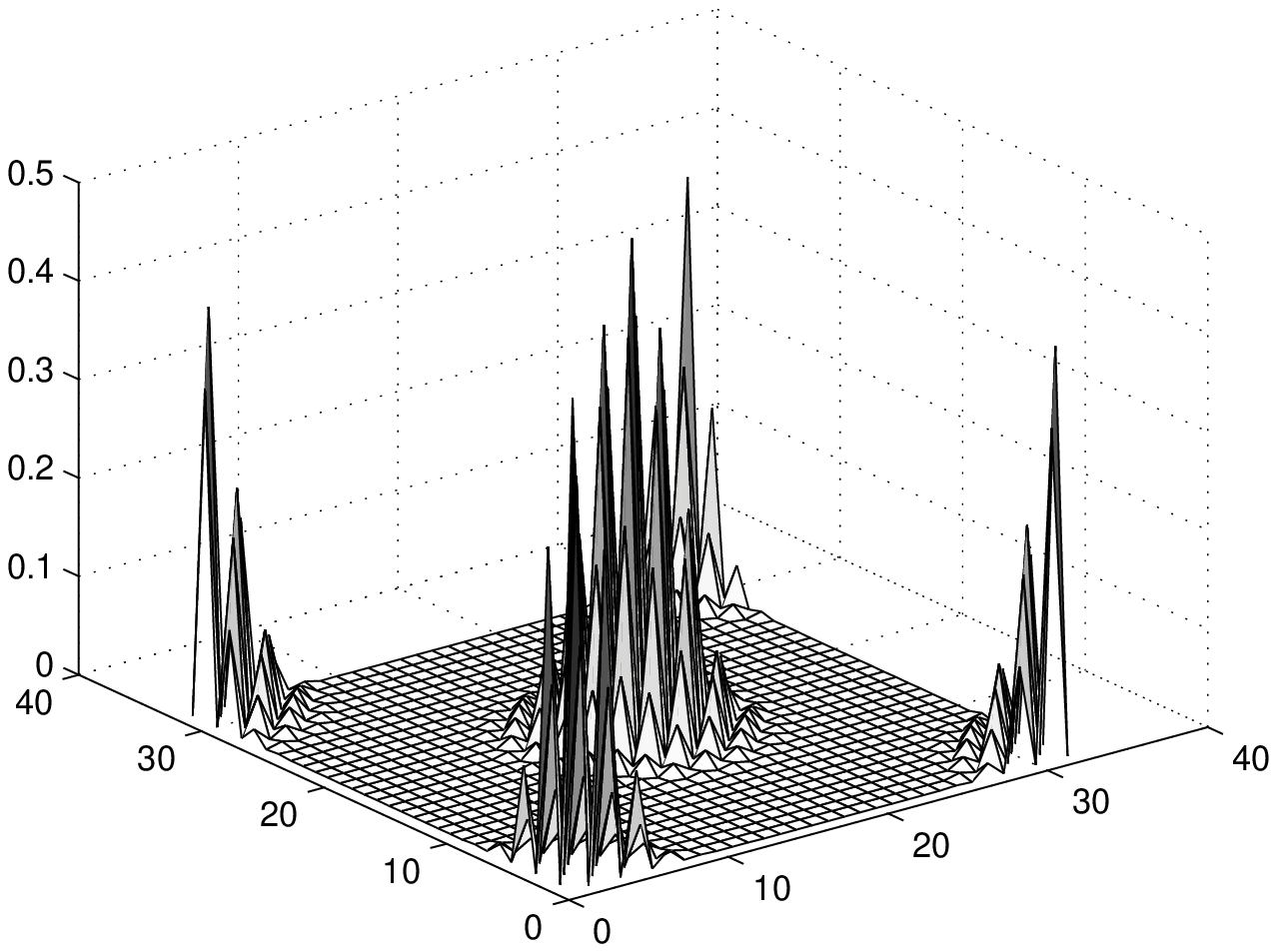}
	\end{center}
	\caption{\small \em Using the algorithm for approximation with frame multipliers in the Gabor case: Approximating the identity by Gabor multiplier with Gauss window ($n=32$) and changing lattice parameters. Top Left: ($a=2$,$b=2$), Top Right: ($a=4$,$b=4$), Bottom Left: ($a=8$,$b=8$), Bottom Right: ($a=16$,$b=16$)} \label{fig:kap1apprgab1}
\end{figure}

\Bsp{sec:testapprframmultgab} We will now use this algorithm for the approximation by regular Gabor multipliers as presented in \cite{feiham1}.
We are using a Gauss window in $\CC^n$ with $n=32$. We are changing the lattice parameters $a$ and $b$. The resulting approximation of the identity can be found in Figure \ref{fig:kap1apprgab1}.

\begin{enumerate}
\item $(g, a = 2, b=2)$. This is nearly a tight frame with the lower frame bound $A = 7.99989$ and the upper frame bound $B = 8.00011$. As expected the identity is approximated very well.
\item $(g,a = 4, b =4)$ : This frame is not tight anymore, as $A = 1.66925$ and $B = 2.36068$ and we can see that the approximation is deviating from identity.
\item $(g,a = 8,b = 8)$ : This is not a frame anymore, but a Bessel sequence with $B = 1.18034$. At least some of the structure (the diagonal dominance) is still kept.
\item $(g,a = 16,b = 16)$ : This is not a frame anymore, but a Bessel sequence with $B = 1.00001$. All structure is (more or less) lost.
\end{enumerate}
This algorithm is not very efficient for the Gabor case as the special structure is not used. For the regular case the algorithm presented in \cite{feiham1} is preferable. 
The MATLAB-codes for these examples can be found  
in the file \verb_testapprGabmult.m_.

\section{Perspectives}

The efficient calculation of the best approximation, in the Hilbert-Schmidt sense, of an arbitrary matrix by a frame multiplier is important for improving the numerical aspects of a frame, e.g.. quotient of the frame bounds, by using weights. Refer to \cite{bogdvan1}, where weighted frames are introduced. This connection is worth to be established more in the future.

The idea for the approximation of Hilbert Schmidt operators can be applied for irregular Gabor multipliers. Even if a lot of the structure is missing, some special properties can be exploited and the algorithms can be sped up. For first ideas refer to \cite{xxlphd1}.


\end{document}